\documentclass[12pt]{article}

\usepackage{amssymb}
\usepackage{amsmath}
\usepackage[latin1]{inputenc}

\newlength{\labst}

\newcommand{\abst}{\vspace{\labst}}

\newenvironment{Proof}{{\bf Proof:}}{\hfill q.e.d.}

\newtheorem{Theorem}{Theorem}

\newenvironment{Example*}{\begin{Example}\DP}{\end{Example}\abst}
\newenvironment{Remark*}{\begin{Remark}\DP}{\end{Remark}\abst}

\DeclareMathAlphabet{\mathsf}{OT1}{cmss}{n}{n}
\SetMathAlphabet{\mathsf}{bold}{OT1}{cmss}{bx}{n}

\begin{document}
\centerline {\Large \bf In which spaces is every curve Lebesgue-Pettis integrable?}
\vskip1.5cm
\centerline{Heinrich v. Weizs\"acker}
\centerline{Kaiserslautern, \today}
\vskip1cm
\centerline{\it Dedicated to the memory of Erik G. F. Thomas} 
\vskip1.5cm
{\bf 1. Introduction.} Around 1980 P. Dierolf und  J. Voigt raised the question which locally convex real vector spaces $E$ are sufficiently 'complete' to have the property in the title. The following text is a  revised translation of an unpublished manuscript of that time.  

Our primary result is that the property in the title (for continuous curves) is equivalent to the 'metric convex compactness property' in the sense of Voigt (1992, \cite{Voigt}): the closed convex hull of a compact metrizable subset is again compact.  There (\cite{Voigt}, Theorem 0.1) an observation of H. Pfister (1981) is reproduced to the effect that this metric convex compactness property is equivalent to the Pettis integrability of all $E$-valued continuous functions on compact metric spaces with a Borel measure $\mu$. So our main point is that it suffices to restrict the attention to the unit interval and the Lebesgue measure. In \cite{Voigt} the metric convex compactness property is compared to other completeness properties like the so-called Mackey-completeness; moreover in that paper it is shown that various spaces of linear operators, e.g. the compact operators in a Banach space, do have this property with respect to the strong operator topology. 

We also show that the property in the title is equivalent to another property which, in view of earlier work of E. G. F. Thomas (1975, e.g. \cite{Thomas}), appears to be the right condition to allow for a satisfactory extension of the Bochner integral beyond the domain of Banach space valued functions. 
\bigskip

{\bf 2. The main result.}

\begin{Theorem} 
A real locally convex Hausdorff space $E$ has the metric convex compactness property if and only if for every continuous map $f: [0,1] \rightarrow E$ the Pettis integral $\int^1_0 f(t)\ dt$ exists in $E$.  
\end{Theorem}

\begin{Proof} We first remark that a Borel map from a topological space $X$ into $E$ is $\nu$-Pettis integrable for a Borel measure $\nu$ on $X$ if and only if the image measure $\mu = \nu \circ f^{-1}$ has a barycenter (resultant) $r(\mu)$ in $E$. In fact $$\langle e', r(\mu) \rangle = \int \langle e', x \rangle \ d\mu = \int_X \langle e', f(t) \rangle\ d\mu(t)$$ for all $e' \in E'$ and thus the Pettis integral of $f$ if it exists is indeed just this barycenter. 
\medskip

1. Assume the metric convex compactness property. If $f : [0,1] \to E$ is continuous then the image set $K = f([0,1])$ is compact und metrizable. Its closed convex hull is compact by assumption. Hence 
the image of the Lebesgue measure $\lambda$ on $[0,1]$ has a barycenter and thus $f$ is $\lambda$-Pettis integrable by the preceding remark. 
\medskip 

2. Conversely assume that every continuous curve $f : [0,1] \to E$
is $\lambda$-Pettis integrable. Combining Pfister's remark (Theorem 0.1 in \cite{Voigt}, see the introduction) and the above it suffices to show that every probability measure $\mu$ on $E$ with compact metrizable support $K$ has a barycenter in $E$. 
\medskip

a. First we assume that $\mu = \sigma \circ g^{-1}$ for a continuous 
map $g : [0,1] \to E$ and a non-atomic probability measure $\sigma$ on $[0,1]$ with $supp\ \sigma = [0,1]$.  Let $F(t) := \sigma ((0,t]).$ 
Then $F$ is strictly increasing and $\sigma$ is the image of $\lambda$ under the continuous map $F^{-1}$. Therefore $\mu = (\lambda \circ (F^{-1})^{-1}) \circ g^{-1} =  \lambda \circ (g \circ F^{-1})^{-1}.$ Thus the barycenter of $\mu$ is just the Pettis integral $\int_0^1 g(x)\ d\sigma(x) = \int_0^1 g(F^{-1}(t))\ dt$ which exists by assumption. 
\medskip

b. Next we drop the support condition on $\sigma$ in the previous argument: Suppose that $\mu = \rho \circ g^{-1}$ for a non-atomic 
probability measure $\rho$ on $[0,1]$. Consider $\sigma = \frac 12(\rho + \lambda)$. Then $supp\ \sigma = [0,1]$ and $\rho = 2\sigma - \lambda$. By the previous step 
$$r(\mu) = \int_0^1 g(x)\ d\rho = 2 \int_0^1 g(x)\ d\sigma - \int_0^1 g(x)\ dx \in E.$$
\medskip

c. Now let $\mu$ be a non-atomic probability measure with compact metrizable support $K \subset E$. Then there is a continuous map from the Cantor set $D$ onto $K$. By linear interpolation in the intervals which constitute the complement of the Cantor set this map can be extended to a continuous map $g:[0,1] \to E$. The measure $\rho$ on $[0,1]$ with distribution function $F(t) = \mu\big(g([0,t])\big)$ is non-atomic with image measure $\mu$. Thus according to the previous step $r(\mu) \in E$.
\medskip

d. Finally let $\mu$ be an arbitrary probability measure with compact metrizable support $K \subset E$. The discrete part of $\mu$ can be written in the form $\mu_d = \sum_{i=1}^\infty a_i \delta_{x_i}$ as an infinite combination of point masses $\delta_{x_i}$. Without loss of generality we may assume that $x_i = 0$ for no $i$. We replace each of these point masses by the uniform distribution on the vector interval
$[\frac 12 x_i, \frac 32 x_i]$. Then the resulting probability measure 
is non-atomic, and it is concentrated on the compact metrizable set $[\frac 12, \frac 32] K$. Thus it has a barycenter and this barycenter clearly is also the barycenter of the original $\mu$. This concludes the proof.
\end{Proof} 
\bigskip

{\bf 3. An Extension}
\medskip

In Theorem 1 a seemingly weak sufficient condition for the metric convex compactness property was given. In Theorem 2 we show that
the metric convex compactness of $E$ property in turn is sufficient for a natural extension of the Bochner integral for $E$-valued functions in the spirit of E.G.F. Thomas \cite{Thomas}. Following Thomas we work with Lusin-measurable functions rather than with measurable functions on abstract measure spaces.
\medskip

We call a positive Borel-measure $\mu$ on a  topological space $X$ {\bf metrically regular}, if 
$$\mu(B) = \sup \{\mu (K): K \subset B,\ K\,  \text{ compact and metrizable} \}$$ holds for all Borel-sets $B.$ 
A function $g$ is called $\mu${\bf -Lusin-measurable}, if for  every Borel-set $B \subset X$ we have $\mu(B) = \sup \{\mu(K) : K \subset B \text{ compact, }g_{|K} \text{ continuous}\}.$ 
\bigskip

 \begin{Theorem} A real locally convex Hausdorff space $E$ has the metric convex compactness property if and only if the following holds: 
 
 Let $X$ be a topological space and let $\nu$ be a metrically regular Borel measure on $X.$ Let $D \subset E$ be a bounded closed absolutely convex set with the associated gauge functional $ \varphi_D : E \rightarrow \mathbb{R}_+$. Let the function $g: X \rightarrow E$ be $\nu$-Lusin-measurable with $\int_X \varphi_D \circ g \ d\nu < \infty.$ Under these conditions $g$ is $\nu$-Pettis-integrable.   
\end{Theorem}

{ \bf Remark.} If a function $g$ satisfies the assumption in the theorem then it is $\nu$-equivalent to a Borel measurable map $\tilde{g}$ from $X$ into the Banach space induced by the gauge functional $\varphi_D$. This function is easily seen to have separable essential range, and therefore is Bochner integrable. So the Pettis-integral actually can be considered as an extended Bochner integral.   

\begin{Proof} Let us remark first that the gauge functional of a bounded closed absolutely convex set $D$ due to Hahn-Banach has
the representation 
 $$\varphi_D (e) = \sup\{ \langle e', e\rangle : e' \in E', \sup_{x \in D} \langle e', x\rangle \le 1\},$$ hence it is lower semi-continuous and in particular Borel measurable. 
 
Assume now that the condition of the theorem holds. Then indeed for  every Borel measure $\nu$ on a compact metric space $X$ every continuous $E$-valued function is Pettis-integrable. As explained above this implies the metric convex compactness property of $E$. 

Conversely assume that the metric compactness property is satisfied and let $X$, $\nu$, $D$ and $g$ be given as indicated. The real valued function defined by $h(x) := \varphi_D (g(x))$ is $\nu$-Lusin-measurable and even $\nu$-integrable. Replacing, if necessary, $d\nu$  by $h d\nu$ and $g$ by $\frac{g}{h} \cdot 1_{\{h>0\}}$, we may and will assume that $\nu(X) < \infty$ and that $g(X) \subset D$. 
The image measure $\mu = \nu \circ g^{-1}$ is finite, metrically regular and concentrated on $D$.

In order to prove the integrability of $g$ it thus suffices to show that each finite, metrically regular measure $\mu$ on a bounded set $D \subset E$ has a barycenter in $E$. If $D$ is compact we know this from the previous discussion. In the general case, let $K_n \subset D, n \in \mathbb{N}$ be disjoint compact metrizable sets such that $\sum_n \mu(K_n) = \mu(E) < \infty$. There is a sequence $(b_n)$ of positive numbers with $b_n \rightarrow 0$ sufficiently slowly such that $\sum^\infty_{n=1} \frac{\mu(K_n)}{b_n} < \infty.$  Consider the measure $\tilde{\mu}$ defined by 
$$\tilde{\mu}(A) = \sum_{n=1}^\infty \frac{\mu(b_n^{-1}A \cap K_n)}{b_n}.$$
This finite metrically regular measure is concentrated on the set $K = \bigcup_n b_nK_n \cup \{0\}$. This set is compact and metrizable; in fact $\bigcup_n K_n$ is bounded and  $b_n \to 0$ and hence the sets $U_m = \bigcup_{n\ge m} b_nK_n \cup \{0\}$ form a base of neighbourhoods of the point $0$ in the relative topology of $K$. Thus $\tilde{\mu}$ has a barycenter and this barycenter is also the barycenter of $\mu$ because for each 
$e' \in E'$ we have 
$$\int e' \ d\tilde{\mu} = \sum_{n=1}^\infty \frac1{b_n}\int_{K_n} \langle e',b_n x \rangle \ d\mu = \sum_{n=1}^\infty \int_{K_n} \langle e', x \rangle \ d\mu = \int e' \ d\mu.$$ 
\end{Proof}
\medskip

{\bf Acknowledgement.} The author is indebted to J. Voigt for quoting the original manuscript in his 1992 paper \cite{Voigt}; moreover to H. Gloeckner, whose curiosity recently was raised by that citation. This renewed interest triggered the present version of this text. Finally my thanks go again to J\"urgen Voigt for carefully reading the final  version.   
\medskip

\bibliographystyle{plain}

\end{document}